\documentclass{amsart}
\usepackage{amsmath}
\usepackage{latexsym,ulem}
\usepackage{amssymb}
\usepackage{graphics}

\input epsf%
\newtheorem{theorem}{Theorem}[section]

\newtheorem{lemma}{Lemma}[section]
\newtheorem{remark}{Remark}[section]

\newtheorem{Definition}{Definition}[section]
\def\qed{{\hfill{\vrule height4pt width3pt depth2pt}}}

   \long\def\comment#1{}

\def\ad#1{\begin{aligned}#1\end{aligned}}  \def\b#1{\mathbf{#1}}
\def\a#1{\begin{align*}#1\end{align*}} \def\an#1{\begin{align}#1\end{align}}
\def\e#1{\begin{equation}#1\end{equation}} \def\d{\operatorname{div}}
\def\p#1{\begin{pmatrix}#1\end{pmatrix}} 
  \numberwithin{equation}{section}
\numberwithin{table}{section}
\numberwithin{figure}{section}

\def\boxit#1{\vbox{\hrule height1pt \hbox{\vrule width1pt\kern1pt
     #1\kern1pt\vrule width1pt}\hrule height1pt }}
 \def\lab#1{\boxit{\small #1}\label{#1}}
  
 \def\meqref#1{\boxit{\small #1}\eqref{#1}}

  \def\lab#1{\label{#1}}  \def\meqref#1{\eqref{#1}}
\begin{document}

\newcommand{\disp}{\displaystyle}
\newcommand{\eps}{\varepsilon}
\newcommand{\To}{\longrightarrow}
\newcommand{\C}{\mathcal{C}}
\newcommand{\K}{\mathcal{K}}
\newcommand{\T}{\mathcal{T}}
\newcommand{\bq}{\begin{equation}}
\newcommand{\eq}{\end{equation}}
\long\def\comments#1{ #1}
\comments{ }
%%%%%%%%%%%%%%%%%%%%%%%%%%%%%%%%%%%%%%%%%%%%%%%%%%%%%%%%%%%%%%%%%%%%%%%%%%

\title[Finite element  approximations]{Finite element  approximations of  symmetric tensors  on simplicial grids in $\mathbb{R}^n$: the higher order case
    }
\date{}

\author {Jun Hu}
\address{LMAM, School of Mathematical Sciences, and  Beijing International
  Center for Mathematical Research,  Peking University,
  Beijing 100871, P. R. China.  hujun@math.pku.edu.cn}

\thanks{The author was supported by  the NSFC Projects 11271035,  91430213 and 11421101.}

\begin{abstract}
The design of mixed finite element methods in linear elasticity with symmetric stress approximations has been a longstanding open problem
until Arnold and Winther designed the first family of mixed finite elements where the discrete stress space is the space
    of $H(\d,\Omega; \mathbb {S})$---$P_{k+1}$ tensors  whose divergence is a $P_{k-1}$ polynomial on each triangle for $k\geq 2$.
    Such a two dimensional  family was extended,  by Arnold, Awanou and Winther,  to a three dimensional family of mixed elements where
     the discrete stress space is the space  of $H(\d,\Omega; \mathbb {S})$---$P_{k+2}$ tensors, whose divergence is a $P_{k-1}$ polynomial on each tetrahedron for $k\geq 2$.  In this paper,  we are able to  construct, in a unified  fashion, mixed finite element methods with symmetric stress approximations on an arbitrary simplex in $\mathbb{R}^n$ for any space dimension.  On the contrary, the discrete stress space here is the space
     of $H(\d,\Omega;\mathbb {S})$---$P_k$ tensors, and the discrete displacement space here is the space of $L^2(\Omega; \mathbb{R}^n)$---$P_{k-1}$ vectors for
      $k\geq n+1$.  These finite element spaces are defined with respect to an arbitrary
simplicial triangulation of the domain, and can be regarded as extensions to any dimension  of those in two and three  dimensions by Hu and Zhang.

   \vskip 15pt


 \vskip 15pt


\end{abstract}
\maketitle

\section{Introduction}
In  the classical  Hellinger-Reissner
mixed formulation of the elasticity equations, the stress is sought in $H(\d,\Omega;\mathbb {S})$
and the displacement in $L^2(\Omega; \mathbb{R}^2)$  for two dimensions and in  $L^2(\Omega; \mathbb{R}^3)$ for three dimensions.
The constructions of stable mixed finite elements using
polynomial shape functions are  a long-standing and challenging problem, see \cite{Arnold2002,Arnold-Awanou-Winther}.
To overcome this  difficulty,  earliest  works adopted composite element
   techniques or weakly symmetric methods,  cf.  \cite{Amara-Thomas,
   Arnold-Brezzi-Douglas, Arnold-Douglas-Gupta,  Johnson-Mercier,
   Morley, Stenberg-1, Stenberg-2, Stenberg-3}.
In \cite{Arnold-Winther-conforming},  Arnold and Winther
   designed the first family of
    mixed finite element methods in 2D,
     based on  polynomial shape  function spaces.
   From then on,   various  stable  mixed elements have been constructed,
        see \cite{Adams-Cockburn,Arnold-Awanou,Arnold-Awanou-Winther,
    Arnold-Winther-conforming,Awanou, Chen-Wang,
     Arnold-Winther-n, Gopalakrishnan-Guzman-n,Hu-Shi,
    Man-Hu-Shi, Yi-3D, Yi, Arnold-Falk-Winther,
    Boffi-Brezzi-Fortin, Cockburn-Gopalakrishnan-Guzman,
    Gopalakrishnan-Guzman, Guzman,Hu-Man-Zhang2014,Hu-Man-Zhang2013}.

For first order systems with symmetric tensors in any space dimension, as the displacement $u$ is in $L^2(\Omega; \mathbb{R}^n)$,
  a natural discretization   is the piecewise $P_{k-1}$ polynomial without interelement continuity.
Even  for two and three dimensional cases,  it is  a  surprisingly hard problem
   if the  stress tensor can be discretized by  an appropriate $P_k$ finite element subspace of
  $H(\d,\Omega;\mathbb {S})$. In fact, in \cite{Arnold-Winther-conforming},   Arnold and Winther
    designed the first family of mixed finite elements where the discrete stress space is the space
    of $H(\d,\Omega; \mathbb {S})$---$P_{k+1}$ tensors  whose divergence is a $P_{k-1}$ polynomial on each triangle with $k\geq 2$; see also \cite{Arnold2002}.
    Such a two dimensional  family was extended to a three dimensional family of mixed elements where
     the discrete stress space is the space  of $H(\d,\Omega; \mathbb {S})$---$P_{k+2}$ tensors with $k\geq 2$; while
       the lowest order element with $k=2$ was first proposed in \cite{Adams-Cockburn}.
In very recent papers \cite{Hu-Zhang2014a} and \cite{Hu-Zhang2014b},  Hu and Zhang  attacked this  open problem by constructing a  suitable $H(\d,\Omega;\mathbb {S})$---$P_k$, instead of above $P_{k+1}$ (2D, $k\geq 3$) or $P_{k+2}$ (3D $k\geq 4$),
    finite element space for the stress discretization.  The analysis there is based on a new idea for analyzing the
        discrete inf--sup condition. More precisely, they first decomposed the discontinuous displacement space into a subspace containing lower order polynomials and its  orthogonal complement space. Second they found that the discrete stress space
        contains  the full $C^0$-$P_k$ space  and  some so-called $H(\d)$  bubble function space on each triangle (2D) or tetrahedron (3D). Third they proved that the full $C^0$-$P_k$ space can control the  subspace containing lower order polynomials while the $H(\d)$  bubble function space is able to deal with  that orthogonal complement space.
         We refer  interested readers to Hu \cite{Hu2013} for similar mixed elements on   rectangular and cuboid  meshes.

  The purpose of this paper is to  generalize, in a unified fashion,  the elements in \cite{Hu-Zhang2014a} and \cite{Hu-Zhang2014b}   to any  dimension.
  In addition, we define a set of local degrees of freedom  for shape function spaces of stress on each element. The analysis here is based on three key ingredients.  First, based on  the tangent vectors of a simplex, we construct $\frac{n(n+1)}{2}$ symmetric matrices of rank one and prove that they are linearly
    independent and  consequently form a basis of the space $\mathbb{S}$.  Second,  by using these matrices of rank one, we
     define a $H(\d)$ bubble function space consisting of polynomials of degree $\leq k$  on each element and
     prove that it is indeed the full $H(\d)$ bubble function space of order $k$.  Third,  we show that   the divergence  space of  the $H(\d)$ bubble function space  is equal to the orthogonal complement space of    the rigid motion space with respect to the  discrete displacement on each element. We stress that  such a result holds     for any $k\geq 2$.

    The rest of the paper is organized as follows.  In the next section, we define finite element spaces
     of symmetric tensors in any space dimension, present a crucial structure of them,
       and a set of local degrees of freedom of shape function spaces on each element.
      We also prove  that  the divergence of  the $H(\d)$ bubble function space  is equal to the orthogonal complement space of    the rigid motion space with respect to the  discrete displacement on each element.
         In Section 3, we apply these spaces  to  first order systems with symmetric tensors
          and prove the well--posedness of the discrete problem.    The paper ends with Section 4, which gives some conclusion.

\section{Finite elements for symmetric tensors}

We consider  mixed finite element methods of
 first order systems with symmetric tensors:
Find $(\sigma,u)\in\Sigma\times V :=H({\rm div},\Omega;
    \mathbb {S})
        \times L^2(\Omega;\mathbb{R}^n)$, such that
\an{\left\{ \ad{
  (A\sigma,\tau)+({\rm div}\tau,u)&= 0 && \hbox{for all \ } \tau\in\Sigma,\\
   ({\rm div}\sigma,v)&= (f,v) &\qquad& \hbox{for all \ } v\in V. }
   \right.\lab{eqn1}
}
Here the symmetric tensor space for the stress $\Sigma$ is  defined by
  \an{   \lab{S}
  H({\rm div},\Omega; \mathbb {S})
    &:= \Big\{ \tau=\p{\tau_{11} &\cdots &\tau_{1n} \\
        \vdots  & \vdots & \vdots \\
        \tau_{n1} &\cdots   & \tau_{nn}  }
     \in H(\d, \Omega;\mathbb{R}^{n\times n}), \ \tau ^T = \tau  \Big\},}
and  the   space for the vector displacement  $V$ is
  \an{  \lab{V}
     L^2(\Omega;\mathbb{R}^n) &:=
     \Big\{ \p{u_1, &\cdots, &u_n}^T, \ u_i \in L^2(\Omega), i=1, \cdots, n \Big\}  .}
	This paper denotes by $H^k(T,X)$ the Sobolev space consisting of
functions with domain $T\subset\mathbb{R}^n$, taking values in the
finite-dimensional vector space $X$, and with all derivatives of
order at most $k$ square-integrable. For our purposes, the range
space $X$ will be either $\mathbb{S},$ $\mathbb{R}^n,$ or
$\mathbb{R}$. Let
$\|\cdot\|_{k,T}$ be the norm of $H^k(T)$. Let $\mathbb{S}$ denote
the space of symmetric tensors, and $H({\rm div},T,\mathbb{S})$
consist of square-integrable symmetric matrix fields with
square-integrable divergence. The H(div) norm is defined by
$$\|\tau\|_{H({\rm div},T)}^2:=\|\tau\|_{0, T}^2
   +\|{\rm div}\tau\|_{0, T}^2.$$
Let $L^2(T,\mathbb{R}^n)$ be the space of vector-valued functions
which are square-integrable.
Here, the compliance tensor
$A=A(x)\in L^\infty(\mathbb{S};\mathbb{S})$, characterizing the
properties of the material, is bounded and symmetric positive
definite uniformly for $x\in\Omega$.

Suppose that the  domain $\Omega$ is subdivided by a family of shape regular simplicial
 grids  $\mathcal{T}_h$ (with the grid size $h$). This paper denotes $P_k(K; X)$ as the space  of polynomials of degree $\leq k$, taking value in the space $X$.

   \subsection{A new basis of the symmetric matrices}
    Let $\b x_0, \cdots, \b x_n$  be the  vertices of
    simplex  $K$.
The referencing mapping is then
\a{  \b x: &= F_K(\hat {\b x})
     = \b x_0 + \p{ \b x_1 -\b x_0, &
                     \cdots,&
                     \b x_n -\b x_0  } \hat{\b x}, }
      mapping the reference tetrahedron $\hat K:=\{ 0 \le \hat x_1, \cdots, \hat x_n,
          1-\sum\limits_{i=1}^n\hat x_i \le 1 \}$ to $K$.
   Then the inverse mapping is
\an{ \lab{imap}  \hat {\b x}: &= \p{\b \nu_1^T\\ \vdots \\ \b \nu_n^T} (\b x-\b x_0), }
where
 \an{\lab{nnt} \p{\b \nu_1^T\\ \vdots \\ \b \nu_n^T }=\p{ \b x_1 -\b x_0, &
      \cdots,  &
                     \b x_n -\b x_0}^{-1}. }
By \meqref{imap}, these normal vectors
        are coefficients of the barycentric variables:
  \a{ \lambda_1(\b x): &=  \b \nu_1 \cdot (\b x-\b x_0), \\
     \vdots & \\
      \lambda_n(\b x): &= \b \nu_n \cdot (\b x-\b x_0), \\
    \lambda_0(\b x): &= 1-\sum\limits_{i=1}^n\lambda_i. }
For any edge $\b x_i\b x_j$ of element $K$, $i\not =j$, let $\b t_{i,j}$ denote  associated tangent vectors, which allow for us to
 introduce the following symmetric matrices of rank one
 \begin{equation}\label{Tnew}
 T_{i,j}:=\b t_{i,j}\b t_{i,j}^T,  0\leq i< j\leq n.
 \end{equation}
For these matrices of rank one, we have the following important result.
 \begin{lemma}\lab{6T}
   The $\frac{(n+1)n}{2}$ symmetric tensors $T_{i,j}$ in \meqref{Tnew} are linearly independent, and
    form a basis of $\mathbb{S}$.
\end{lemma}

\begin{proof}  Each matrix $T_{i,j}=\b t_{i,j}  \b t_{i,j}^T$ is a positive
  semi-definite matrix, on a simplex $K$.
   We would show that the constants $c_{i,j}$ are all
      equal to zero in
  \a{  \tau= \sum\limits_{0\leq i<j\leq n}c_{i,j} T_{i,j}=0. }
  Let  $\b \nu_0$ be the normal  vector to the $n-1$ dimensional simplex $\triangle_{n-1}\b x_1\cdots \b x_n$. This leads to
    \begin{equation}
  \b \nu_0^T \b t_{i,j}=0, 1\leq i<j\leq n,
  \end{equation}
  and
  \begin{equation}\label{eq2.6}
 \b \nu_0^T\b t_{0, j}\not =0,  1\leq j \leq  n.
 \end{equation}
  This gives
 \begin{equation}
 \begin{split}
\b \nu_0^T \tau &=\b \nu_0^T\sum\limits_{0\leq i<j\leq n}c_{i,j}\b t_{i, j}\b t_{i,j}^T
=\b \nu_0^T\sum\limits_{1\leq j\leq n}c_{0,j}\b t_{0, j}\b t_{0,j}^T\\
&=\sum\limits_{1\leq j\leq n}\tilde{c}_{0,j}\b t_{0,j}^T=0,
\end{split}
 \end{equation}
 where  $\tilde{c}_{0,j}=c_{0, j}\b \nu_0^T\b t_{0, j}$.    Since $\b t_{0, j}$, $1\leq j\leq n$, are linearly independent,  this yields
  \begin{equation}\label{c1j}
 \tilde{c}_{0,j}=0, 1\leq j\leq n.
  \end{equation}
    This and \eqref{eq2.6} yield
 \begin{equation}
c_{0,j}=0,1\leq j\leq n.
 \end{equation}
 A similar argument by using $\b \nu_i$, $i\not =0$, proves the desired result.
\end{proof}

\subsection{The bubble--function space }
With these symmetric matrices $T_{i,j}$ of rank one,  we define a $H(\d, K; {\mathbb S})$ bubble function space
\an{\lab{bh}  \Sigma_{K, k,  b}: = \sum\limits_{0\leq i<j\leq n} \lambda_i\lambda_jP_{k-2}(K; \mathbb{R})T_{i, j}.}
Then we  propose to define the  discrete stress space $\Sigma_{k,h}$ which has the following {\it crucial} structure:
 \an{\lab{Sh}
 \Sigma_{k, h}: =\Big\{~\sigma
    &\in H(\d,\Omega,\mathbb{S}),  \sigma=\sigma_c+\sigma_b,
             \ \sigma_c\in H^1(\Omega, \mathbb{S}), \\
     & \ \sigma_c|_K\in P_k(K, \mathbb{S})
	        \,,
	\nonumber  \ \sigma_b|_K\in\Sigma_{K, k, b},  \forall K\in\mathcal{T}_h \Big\}, }
 which is a $H(\d )$ bubble enrichment of the $H^1$ space $\widetilde \Sigma_{k,h}:=\Sigma_{k,h}\cap H^1(\Omega; \mathbb{S})$ of $\Sigma_{k,h}$.
 This generalizes  the results of \cite{Hu-Zhang2014a, Hu-Zhang2014b} for  both two and three dimensions to  the general case in any space dimension.
 Such a structure  has already enabled  us to write down directly  the basis of  $\Sigma_{k, h}$;  see  \cite{Hu-Zhang2014a,Hu-Zhang2014b} for more details in both two and three dimensions.  Next we  plan,  as  it has been done for most of usual finite element methods in the literature,
  to  define a set of  local degrees of freedom of shape function spaces $P_k(K, \mathbb{S})$ on each element.
  To this end, we define the full $H(\d, K; {\mathbb S})$ bubble function space consisting of polynomials of degree $\leq k$
\begin{equation}
\Sigma_{\partial K,  k, 0}:=\{\tau\in P_{k}(K; \mathbb{S}), \tau\b \nu|_{\partial K}=0 \}.
\end{equation}
Here $\nu$ is the normal vector of $\partial K$.

   \begin{lemma}\label{equiv}  It holds that
   \begin{equation}
   \Sigma_{K, k, b}=\Sigma_{\partial K,  k, 0}.
   \end{equation}
      \end{lemma}
\begin{proof} Consider a function  $\tau\in\lambda_i\lambda_jP_{k-2}(K; \mathbb{R})T_{i, j}$, $0\leq i<j\leq n$.
Note that $\tau$  vanishes on the $n-1$ dimensional  simplices
$$\triangle_{n-1}\b x_0\cdots \b x_{i-1}\b x_{i+1}\cdots\b x_j\cdots \b x_n,$$
 $$\triangle_{n-1}\b x_0\cdots \b x_i\cdots\b x_{j-1}\b x_{j+1}\cdots \b x_n.$$
For any $n-1$ dimensional simplex which takes edge $\b x_i \b x_j$, its normal vector, say $\b \nu$, is perpendicular to the tangent vector $\b t_{i,j}$ of  edge $\b x_i \b x_j$, which implies that $\tau \b \nu=0$ on such a $n-1$ dimensional simplex and consequently $\tau \in \Sigma_{\partial K,  k, 0}$.
 Hence
 \begin{equation}\label{in}
 \Sigma_{K, k, b}\subset \Sigma_{\partial K,  k, 0}.
  \end{equation}
  Next we  show the converse of \eqref{in}.  Given $\tau\in \Sigma_{\partial K,  k, 0}$,  the boundary condition $\tau\b \nu|_{\partial K}=0$ indicates that $\tau$ vanishes at all the vertices of  $K$.  Let $\mathbb{N}_b$ denote all the nodes
    except the vertices of $K$ for the space $P_k(K; \mathbb{R})$.  Given $\mathbb{P}_\ell\in \mathbb{N}_b$, let $\varphi_\ell\in P_k(K; \mathbb{R})$ denote the usual   associated  nodal  Lagrange basis function, namely, $\varphi_\ell(\mathbb{P}_\ell)=1$ and
    $\varphi_\ell$  vanishes at all the  other nodes for the space  $P_k(K; \mathbb{R})$. It follows from Lemma \ref{6T} that
  \begin{equation}\label{14}
  \tau=\sum\limits_{0\leq i< j\leq n} T_{i,j}\bigg(\sum\limits_{\mathbb{P}_\ell\in \mathbb{N}_b}c_{\ell,i,j}\varphi_\ell\bigg).
  \end{equation}
 Note that $\varphi_\ell$ has a homogeneous  expression by $\lambda_0, \cdots, \lambda_n$.  Therefore,  we have
 \begin{equation}\label{15}
 \sum\limits_{\mathbb{P}_\ell\in \mathbb{N}_b}c_{\ell,i,j}\varphi_\ell=\sum\limits_{m_0+m_1+\cdots+m_n= k}c_{(ij), m_0, m_1, \cdots, m_n}\lambda_0^{m_0}\cdots\lambda_n^{m_n}.
 \end{equation}
 We claim that $\sum\limits_{\mathbb{P}_\ell\in \mathbb{N}_b}c_{\ell,i,j}\varphi_\ell$ has a  factor
 $\lambda_i\lambda_j$, namely,
 \begin{equation}
 \sum\limits_{\mathbb{P}_\ell\in \mathbb{N}_b}c_{\ell,i,j}\varphi_\ell=\lambda_i\lambda_j\sum\limits_{m_0^\prime+m_1^\prime+\cdots+m_n^\prime= k-2}c_{(ij), m_0\prime, m_1^\prime, \cdots, m_n^\prime}^\prime\lambda_0^{m_0^\prime}\cdots\lambda_n^{m_n^\prime}.
 \end{equation}
 Without loss of generality, we consider the case where $i=0$ and $j=1$. Suppose that there is a term $f_1T_{0,1}$ such that $f_1$ is a polynomial of degree $\leq k$ and does not contain a factor $\lambda_0$.  Next we shall show that $f_1=0$.  In fact,  all the terms of \eqref{14} which do not contain the  factor $\lambda_0$ and whose normal  components (namely $T_{0,j}\nu_0\not=0$, $\nu_0$ is the normal vector of $\triangle_{n-1}\b x_1\cdots\b x_n$) do not vanish on the $n-1$ dimensional simplex
  $\triangle_{n-1}\b x_1\cdots\b x_n$  can be expressed as
 \begin{equation}
 \sum\limits_{j=1}^n f_j T_{0,j},
 \end{equation}
 where $f_j$, $j=1, \cdots, n$, are polynomials of degree $\leq k$. Since $f_j$ do not contain the factor $\lambda_0$,
  it is of the form
  \begin{equation}\label{expression}
  f_j=\sum\limits_{r_1+\cdots +r_n= k}c_{j, r_1,\cdots, r_n}\lambda_1^{r_1}\cdots \lambda_n^{r_n}.
  \end{equation}
  Since $\tau \b \nu_0=0$ on the $n-1$ dimensional simplex
  $\triangle_{n-1}\b x_1\cdots\b x_n$,
  \begin{equation}
  \sum\limits_{j=1}^n (\b t_{0, j}^T \b \nu_0)f_j \b t_{0,j}\bigg|_{\triangle_{n-1}\b x_1\cdots\b x_n}=0.
  \end{equation}
  Since,  for $j=1, \cdots, n$,  $\b t_{0,j}^T\b \nu_0\not=0$,   and $\b t_{0, j}$ are linearly independent,   this leads to
   \begin{equation}
   f_j|_{\triangle_{n-1}\b x_1\cdots \b x_n}\equiv 0.
   \end{equation}
   Note that $\lambda_1^{r_1}\cdots \lambda_n^{r_n}|_{\triangle_{n-1}\b x_1\cdots \b x_n}$, $\sum\limits_{i=1}^nr_i= k$,
   form a basis of $P_k(\triangle_{n-1}\b x_1\cdots\b x_n; \mathbb{R})$. This and the above equation show that
   \begin{equation}
   c_{j, r_1,\cdots, r_n}=0.
   \end{equation}
     This, in turn, implies that
     \begin{equation}
      f_j\equiv 0.
     \end{equation}
   Therefore $f_1=0$ which implies that all the terms on the right hand side of  \eqref{15}  has a factor $\lambda_0$.
    A  similar argument shows that   all the terms on the right hand side of  \eqref{15}  has a factor $\lambda_1$.
   Hence
  \begin{equation}
  \tau\in \Sigma_{K, k, b}.
  \end{equation}
  This completes the proof.
\end{proof}

\subsection{Degrees of freedom}
Before we define the degrees of freedom, we need a classical result and its variant.
\begin{lemma} It holds the following  Chu-Vandermonde combinatorial identity and its variant
\begin{equation}\label{van}
\sum\limits_{\ell=0}^nC_{n+1}^{\ell+1}C_{k-1}^\ell =\sum\limits_{\ell=0}^nC_{n+1}^{n-\ell}C_{k-1}^\ell=C_{n+k}^n,
\end{equation}
\begin{equation}\label{varvan}
\sum\limits_{\ell=0}^nC_{n+1}^{\ell+1}C_{k-1}^\ell C_{\ell+1}^2=\frac{(n+1)n}{2}C_{n+k-2}^n,
\end{equation}
where the combinatorial number $C_{n}^m=\frac{n\cdots (n-m+1)}{m\cdots 1}$ for $n\geq m$ and $C_n^m=0$ for $n<m$.
\end{lemma}
\begin{proof}  The identity \eqref{van} is  classical, and \eqref{varvan} is its  variant.  For readers' convenience,
  we sketch the proof for them. It follows from the well-known binomial theorem that
  \begin{equation*}
  (1+t)^{n+1}(1+t)^{k-1}=(1+t)^{n+k}=\sum\limits_{m=0}^{n+k}C_{n+k}^mt^m.
  \end{equation*}
  On the other hand, we have
  \begin{equation*}
  (1+t)^{n+1}(1+t)^{k-1}=\sum\limits_{m_1=0}^{n+1}C_{n+1}^{m_1}t^{m_1}
  \sum\limits_{m_2=0}^{k-1}C_{k-1}^{m_2}t^{m_2}.
    \end{equation*}
  A combination of these two  equations leads to
  \begin{equation*}
  \begin{split}
  C_{n+k}^n&=\sum\limits_{m_1+m_2=n}C_{n+1}^{m_1}C_{k-1}^{m_2}
  =\sum\limits_{\ell=0}^nC_{n+1}^{n-\ell}C_{k-1}^{\ell},
  \end{split}
  \end{equation*}
  which proves \eqref{van}. To prove \eqref{varvan}, we consider
  \begin{equation*}
\sum\limits_{\ell=0}^nC_{n+1}^{\ell+1}C_{k-1}^\ell \frac{(\ell+1)\ell}{(n+1)n}=\sum\limits_{\ell=0}^nC_{n-1}^{\ell-1}C_{k-1}^\ell
=\sum\limits_{\ell=0}^nC_{n-1}^{n-\ell}C_{k-1}^\ell=C_{n+k-2}^n.
\end{equation*}
\end{proof}

\begin{theorem}\label{main1} A  matrix field  $\tau\in P_k(K; \mathbb{S})$ can be uniquely determined by the  degrees of
 freedom  from (1) and (2)
\begin{enumerate}
 \item For each $\ell$ dimensional simplex $\triangle_\ell$ of $K$, $0\leq \ell\leq n-1$,  with $\ell$ linearly independent tangential vectors $\b t_1, \cdots, \b t_\ell$, and $n-\ell$ linearly  independent normal vectors $\b \nu_1, \cdots, \b \nu_{n-\ell}$, the mean moments of degree at most $k-\ell-1$ over $\triangle_\ell$, of~~ $\b t_l^T\tau \b \nu_i$, $\b \nu_i^T\tau\b \nu_j$, $l=1, \cdots, \ell$, $i,j=1, \cdots, n-\ell$,  $\big( C_{n+1-\ell}^2+\ell (n-\ell)\big) C_{k-1}^\ell=\frac{(n-\ell)(n+\ell+1)}{2}C_{k-1}^\ell$ degrees of freedom  for each $\triangle_\ell$;
 \item the values $\int_K \tau: \theta d\b x$ for any $\theta \in \Sigma_{K, k, b}$, $\frac{(n+1)n}{2}C_{n+k-2}^n$ degrees of freedom.
 \end{enumerate}
  \end{theorem}
\begin{proof} We assume that all degrees of freedom vanish and show that $\tau=0$. Note that the mean moment becomes the value of $\tau$ for a $0$ dimensional simplex $\triangle_0$, namely, a vertex, of $K$.  The  first set of degrees of freedom imply that $\tau \b \nu =0$ on $\partial K$. Then the second set of degrees of freedom and Lemma \ref{equiv} show $\tau=0$. Since the
  number of degrees of freedom in the second set follows immediately from  Lemma \ref{6T}, we only need to  prove
  the number of degrees of freedom in the first set.  The number of  $\b t_l^T\tau \b \nu_i$, $l=1, \cdots, \ell$, $i=1, \cdots, n-\ell$, is
  $$
  \ell(n-\ell),
  $$
  while,  by symmetry, the number of  $\b \nu_i^T\tau\b \nu_j$, $i,j=1, \cdots, n-\ell$, reads
  $$
  \frac{(n-\ell)(n-\ell+1)}{2}.
  $$
  The number of the mean moments of degree at most $k-\ell-1$ over $\triangle_\ell$ is $C_{k-1}^\ell$. These  imply the
   number of degrees of freedom in the first set is
   $$
   (\ell(n-\ell)+ \frac{(n-\ell)(n-\ell+1)}{2})C_{k-1}^\ell=\frac{(n-\ell)(n+\ell+1)}{2}C_{k-1}^\ell.
   $$
 Hence   the sum of  degrees of freedom in both sets  reads
  \begin{equation*}
  \begin{split}
  &\sum\limits_{\ell=0}^{n-1}C_{n+1}^{\ell+1}\frac{(n-\ell)(n+\ell+1)}{2}C_{k-1}^\ell+\frac{(n+1)n}{2}C_{n+k-2}^n\\
  &=\frac{(n+1)n}{2}\sum\limits_{\ell=0}^{n-1}C_{n+1}^{\ell+1}C_{k-1}^\ell-\sum\limits_{\ell=0}^{n-1}\frac{\ell(\ell+1)}{2}C_{n+1}^{\ell+1}C_{k-1}^\ell
  +\frac{(n+1)n}{2}C_{n+k-2}^n\\
  &=\frac{(n+1)n}{2}\sum\limits_{\ell=0}^{n}C_{n+1}^{\ell+1}C_{k-1}^\ell-\sum\limits_{\ell=0}^{n}\frac{\ell(\ell+1)}{2}C_{n+1}^{\ell+1}C_{k-1}^\ell
  +\frac{(n+1)n}{2}C_{n+k-2}^n.
  \end{split}
  \end{equation*}
  Then it  follows from the Chu-Vandermonde combinatorial identity \eqref{van} and its variant \eqref{varvan} that it is equal to
  $\frac{n(n+1)}{2}C_{n+k}^n$ the dimension of $P_k(K;  \mathbb{S})$.
\end{proof}

\begin{remark} It follows from Theorem \ref{main1} that, for any dimension, if $k=1$,
$\Sigma_{k,h}$ becomes a  $H^1$ conforming  approximation of $\Sigma:=H(\d, \Omega; \mathbb{S})$.
For one dimensional case with $n=1$, for any $k$, $\Sigma_{k,h}$  becomes the usual $H^1$ finite element space of degree $k$.
\end{remark}

\subsection{The divergence space of the bubble function space}
Before ending this section,  we prove an important result concerning the divergence space of the bubble function space.
 To this end, we introduce the following rigid motion space on each element $K$.
  \begin{equation}\label{R}
  R(K):=\{ v\in H^1(K; \mathbb{R}^n), (\nabla v+\nabla v^T)/2=0\}.
  \end{equation}
It follows from the definition that  $R(K)$ is a subspace of $P_1(K; \mathbb{R}^n)$.
For $n=1$, $R(K)$ is the constant function space over $K$. The dimension of $R(K)$ is $\frac{n(n+1)}{2}$. This
 allows for defining the orthogonal complement  space of  $R(K)$ with  respect to $P_{k-1}(K; \mathbb{R}^n)$ by
 \begin{equation}
 R^\perp(K):=\{v\in P_{k-1}(K; \mathbb{R}^n), (v, w)_K=0\text{ for any }w\in R(K)\},
 \end{equation}
 where the inner product $(v, w)_K$ over $K$ reads $(v, w)_K=\int_K v\cdot w d\b x$.
\begin{theorem}\label{equal} For any $K\in\mathcal{T}_h$, it holds that
\begin{equation}
\d \Sigma_{K, k, b}=R^\perp(K).
\end{equation}
\end{theorem}
\begin{proof} For any $\tau\in \Sigma_{K, k, b}$, an integration by parts yields
\begin{equation*}
\int_K\d \tau\cdot w d\b x=0 \text{ for any } w\in R(K).
\end{equation*}
This implies that
\begin{equation}
\d \Sigma_{K, k, b}\subset R^\perp(K).
\end{equation}
Next we show the converse.  In fact,
if $\d \Sigma_{K,k,b} \ne R^\perp(K)$,
    there is a nonzero $ v\in R^\perp (K)$
    such that
  \a{ \int_K \d \tau \cdot v\, d\b x=0 \quad \forall \tau \in \Sigma_{K,k,b}. }
By integration by parts,   for  $\tau \in \Sigma_{K, k, b}$,
  we have
  \an{ \lab{o5}
   \int_K  \d \tau \cdot
     v d\b x = -\int_{ K} \tau: \epsilon(v) d\b x=0, }
  where $\epsilon(v)$ is the symmetric gradient, $(\nabla v + \nabla^T v)/2$.

  By Lemma \ref{6T}, $T_{i,j}$, $0\leq i<j\leq n$   defined in \eqref{Tnew},
 form a basis of the space of symmetric matrices in $\mathbb{R}^{n\times n}$. Then there exists
 an associated  dual  basis, say $M_{i,j}$, $0\leq i<j\leq n$, such that
  \begin{equation}
  T_{i,j}: M_{k,l}=\delta_{i,k}\delta_{j,l}, 0\leq i< j\leq n, 0\leq k<l\leq n.
  \end{equation}
  Here the inner product of two matrices $A=(a_{ij})_{i,j=1}^n$ and $B=(b_{ij})_{i,j=1}^n$ is defined as
  $$
  A: B=\sum\limits_{i=1}^n\sum\limits_{j=1}^na_{ij}b_{ij}.
  $$

As $\epsilon(v)\in  P_{k-2}(K; \mathbb{S})$, it follows  that
there exist $q_{i,j}\in P_{k-2}(K;  \mathbb{R})$, $0\leq i< j\leq n$, such that
  \an{\lab{expan} \epsilon(v)=\sum\limits_{0\leq i<j\leq n} q_{i,j}M_{i, j}. }
Selecting $\tau = \sum\limits_{0\leq i<j\leq n}\lambda_i\lambda_j q_{i,j}T_{i,j}\in \Sigma_{K, k, b}$,
         we have,
   \a{ 0 = \int_K \tau: \epsilon(v) d\b x
         = \sum\limits_{0\leq i<j\leq n}\int_K \lambda_i\lambda_j q_{i,j}^2(\b x) d\b x. }
As $\lambda_i\lambda_j>0$ on $K$,  we conclude that $q_{i, j}\equiv 0$, which implies that $v$ is a rigid motion.
This contradicts with  $v\in R^\perp(K)$. Hence $R^\perp(K)\subset \d \Sigma_{K, k, b}$, which completes the proof.
\end{proof}

\section{Mixed  methods of first order systems with symmetric tensors}\label{sec3}
\subsection{Mixed methods}
We propose to use the space $\Sigma_{k,h}$, with $k\geq n+1$,  defined in \eqref{Sh} to approximate $\Sigma$.  In order  get a stable pair of spaces, we take the discrete displacement space  as the full $C^{-1}$-$P_{k-1}$ space
 \an{ \lab{Vh}
   V_{k,h}: = \{v\in L^2(\Omega; \mathbb{R}^n),
         v|_K\in P_{k-1}(K;  \mathbb{R}^n)\ \hbox{ for all } K\in\mathcal{T}_h \}.
    }

 It follows from the definition of $V_{k,h}$ ($P_{k-1}$ polynomials)
    and $\Sigma_{k,h}$ ($P_k$ polynomials) that
   \a{ \d  \Sigma_{k,h} \subset V_{k,h}.}
This, in turn, leads to a strong divergence-free space:
 \an{ \lab {kernel}
    Z_h&:= \{\tau_h\in\Sigma_{k,h} \ | \ (\d\tau_h, v)=0 \quad
	\hbox{for all } v\in V_{k,h}\}\\
    \nonumber
          &= \{\tau_h \in\Sigma_{k,h} \ | \  \d \tau_h=0
    \hbox{\ pointwise } \}.
    }

The mixed finite element approximation of Problem (1.1) reads: Find
   $(\sigma_h,~u_h)\in\Sigma_{k,h}\times V_{k,h}$ such that
 \e{ \left\{ \ad{
    (A\sigma_h, \tau)+({\rm div}\tau, u_h)&= 0 &&
              \hbox{for all \ } \tau \in\Sigma_{k,h},\\
     (\d\sigma_h, v)& = (f, v) &&  \hbox{for all \ } v\in V_{k,h}.
      } \right. \lab{DP}
    }

\subsection{Stability analysis and error estimates}
The convergence of the finite element solutions follows
   the stability and the standard approximation property.
So we consider first the well-posedness  of the discrete problem
    \meqref{DP}.
By the standard theory,  we only need to prove
   the following two conditions, based on their counterpart at
    the continuous level.

\begin{enumerate}
\item K-ellipticity. There exists a constant $C>0$, independent of the
   meshsize $h$ such that
    \an{ \lab{below} (A\tau, \tau)\geq C\|\tau\|_{H(\d)}^2\quad
       \hbox{for all } \tau \in Z_h, }
    where $Z_h$ is the divergence-free space defined in \meqref{kernel}.

\item  Discrete B-B condition.
    There exists a positive constant $C>0$
            independent of the meshsize $h$, such that
    \an{\lab{inf-sup}
   \inf_{0\neq v\in V_{k,h}}   \sup_{0\neq\tau\in\Sigma_{k,h}}\frac{({\rm
        div}\tau, v)}{\|\tau\|_{H(\d)}  \|v\|_{0} }\geq
    C .}
\end{enumerate}

It follows from $\d  \Sigma_{k,h} \subset V_{k,h}$ that $\d  \tau=0$ for
   any $\tau\in Z_h$. This implies the above K-ellipticity condition
	\meqref{below}.
It remains to show the discrete B-B condition \meqref{inf-sup},
  in the following two lemmas.

  We recall  the subspace $\widetilde \Sigma_{k,h}:=\Sigma_{k,h}\cap H^1(\Omega; \mathbb{S})$ of $\Sigma_{k,h}$. For $\tau\in \widetilde \Sigma_{k,h}$, the degrees of freedom on any element $K$
   are:  for each $\ell$ dimensional simplex $\triangle_\ell$ of $K$, $0\leq \ell\leq n$,  the mean moments of degree at most $k-\ell-1$ over $\triangle_\ell$, of $\tau$.  A standard argument is able to prove that these degrees of freedom are unisolvent.

\begin{lemma}\label{lemma1}
For any $v_h\in V_{k,h}$,  there is a $\tau_h \in
    \widetilde \Sigma_{k,h}$  such that,
  for all polynomial $p\in R(K)$, $K\in\mathcal{T}_h$,
   \bq\label{l-1}  \int_K (\d\tau_h-v_h) \cdot p\, d\b x=0
      \quad \hbox{\rm
      and } \quad \|\tau_h\|_{H(\d)}\leq C\|v_h\|_{0}. \eq
     \end{lemma}
\begin{proof}
  Let $v_h\in V_{k,h}$.
  By the stability of the continuous formulation,
      cf. \cite{Arnold-Winther-conforming} for two dimensional case, there is
    a $\tau \in  H^1(\Omega; \mathbb{S})$ such that,
   \a{ \d\tau=v_h \quad \hbox{\rm
      and } \quad \|\tau\|_{1}\le  C\|v_h\|_{0}. }
   In this paper, we only consider the domain such that the above stability holds. We refer interested authors
    to \cite{GiraultRaviart1986} for the  classical result which states  it  is true for Lipschitz domains in
     $\mathbb{R}^n$; see \cite{Duran2001} for more refined results.
First let $I_h$ be a  Scott-Zhang \cite{Scott-Zhang}
   interpolation operator such that
  \begin{equation}\label{eq3.9}
  \|\tau-I_h\tau\|_0+h\|\nabla I_h\tau\|_0\leq Ch\|\nabla\tau\|_0.
  \end{equation}
  Since $k\geq n+1$, $k-(n-1)-1\geq 1$,  for each $n-1$ dimensional simplex $\triangle_{n-1}$ of $K$,
   there are at least  $n$ bubble functions on $\triangle_{n-1}$ for
    each component of $\tau$.   In fact let $\mathbb{T}_{ij}$, $1\leq i< j\leq n$ be the canonical basis of the space
     $\mathbb{S}$. There are $C_{k-1}^{k-n}$ Lagrange basis functions $\varphi_\ell\in
     \{p\in H^1(\Omega; \mathbb{R}), p|_K\in P_k(K; \mathbb{R}), \text{for any }K\in\mathcal{T}_h\}$, $\ell=1, \cdots, C_{k-1}^{k-n}$, such that $\varphi_\ell$ vanish on $\partial (K^+\cup K^-)$, where $K^+$ and $K^-$ are two  elements
      that share the common $n-1$ dimensional simplex $\triangle_{n-1}$. Then $\varphi_\ell\mathbb{T}_{ij}$,
      $1\leq i< j\leq n$, $\ell=1, \cdots, C_{k-1}^{k-n}$, are matrix--valued  bubble functions,which are linearly independent.
      These bubble  functions allow for defining a correction $\delta_h\in \widetilde \Sigma_{k,h}$ such that
    \begin{equation}
    \int_{\triangle_{n-1}} \delta_h\nu \cdot p d \b s=\int_{\triangle_{n-1}} (\tau-I_h\tau)\nu\cdot p d\b s\text{ for any }
    p\in R(K)|_{\triangle_{n-1}}.
    \end{equation}
        Finally we take
    \begin{equation}
    \tau_h=I_h\tau+\delta_h.
    \end{equation}
    We get  a partial-divergence matching property of $\tau_h$:
     for any $p\in R(K)$,  as the symmetric gradient $\epsilon( p)=0$,
 \a{  \int_K (\d\tau_h-v_h) \cdot p \, d\b x
     & = \int_K (\d\tau_h-\d \tau) \cdot p \, d\b x\\
        & = \int_{\partial K}  (\tau_h-\tau)\nu\cdot p\,  d\b s =0.}
 It remain to show the  stability estimate. It is standard to use a scaling argument and the trace theory to show  that
  \begin{equation*}
  \|\delta_h\|_0+h\|\nabla\delta_h\|_0\leq C\big(\|\tau-I_h\tau\|_0+h\|\nabla(\tau-I_h\tau)\|_0\big).
  \end{equation*}
  Then the stability  estimate in \eqref{l-1} follows from \eqref{eq3.9} and
     the triangle inequality.
  \end{proof}

 We are in the position to show the well-posedness of the discrete problem.
\begin{theorem}
 For the discrete problem (\ref{DP}), the K-ellipticity \meqref{below}
    and the discrete B-B
 condition \meqref{inf-sup} hold uniformly.
  Consequently,  the discrete
     mixed problem \meqref{DP} has a unique solution
         $(\sigma_h,~u_h)\in\Sigma_{k,h}\times V_{k,h}$.
\end{theorem}
\begin{proof}  The  K-ellipticity immediately follows from the fact
      that $\d  \Sigma_{k,h} \subset V_{k,h}$.
      To prove the  discrete B-B  condition \meqref{inf-sup},
      for any $v_h\in V_{k,h}$,
        it follows from Lemma \ref{lemma1} that there exists a
      $\tau_{1}\in \Sigma_{k,h}$ such that,  for any polynomial $p\in R(K)$,
   \bq  \int_K (\d\tau_1-v_h) \cdot pd\b x=0
      \quad \hbox{\rm
      and } \quad \|\tau_1\|_{H(\d)}\leq C\|v_h\|_{0}. \eq
Then it follows from Theorem \ref{equal} that
 there is a $\tau_2\in\Sigma_{k, h}$ such that  $\tau_2|_K \in \Sigma_{K, k, b}$  and
   \begin{equation}
     \d\tau_2 = v_h-\d\tau_1, \|\tau_2\|_0=\min\{\|\tau\|_0, \d\tau=v_h-\d\tau_1, \tau\in \Sigma_{K, k, b} \}
     \end{equation}
     It follows from the definition of $\tau_2$ that $\|\d \tau_2\|_0$  defines a norm for it. Then, a scaling argument proves
\begin{equation}
 \|\tau_2 \|_{H(\d)}\leq C\|\d\tau_1-v_h\|_{0}.
\end{equation}
Let $\tau=\tau_1+\tau_2$.  This implies that
\begin{equation}
\d\tau=v_h \text{ and } \|\tau\|_{H(\d)}\leq C\|v_h\|_{0},
\end{equation}
this proves the discrete B-B condition \meqref{inf-sup}.
\end{proof}

\begin{theorem}\label{MainError} Let
  $(\sigma, u)\in\Sigma\times V$ be the exact solution of
   problem \meqref{eqn1} and $(\tau_h, u_h)\in\Sigma_{k,h}\times
   V_{k,h}$ the finite element solution of \meqref{DP}.  Then,
   for $k\ge n+1$,
\an{ \lab{t1} \|\sigma-\sigma_h\|_{H({\rm div})}
    + \|u-u_h\|_{0}&\le     Ch^k(\|\sigma\|_{k+1}+\|u\|_{k}).
      }
\end{theorem}

\begin{proof}
 The stability of the elements and the standard theory of mixed
  finite element methods \cite{Brezzi, Brezzi-Fortin} give the
  following quasioptimal error estimate immediately
\an{
  \label{theorem-err1} \|\sigma-\sigma_h\|_{H({\rm
  div})}+\|u-u_h\|_{0}\leq C \inf\limits_{\tau_h\in\Sigma_{k,h},v_h\in
  V_{k,h}}\left(\|\sigma-\tau_h\|_{H({\rm div})}+\|u-v_h\|_{0}\right).}
Let $P_h$ denote the local $L^2$ projection operator,
   or  element-wise interpolation operator,  from $V$ to $V_{k,h}$,
  satisfying the error estimate
\an{\label{proj-error}
   \|v-P_hv\|_{0}\leq Ch^k\|v\|_{k} \text{ for any }v\in H^k(\Omega; \mathbb{R}^n). }
Choosing $\tau_h=I_h\sigma\in \Sigma_{k,h}$
    where $I_h$ is defined in \eqref{eq3.9} as $I_h$ preserves symmetric $P_k$ functions locally,
   \an{ \lab{p-err2}
       \|\sigma -\tau_h\|_{0} + h |\sigma -\tau_h|_{H(\d)}
        \le Ch^{k+1} \|\sigma\|_{k+1}. }
 Let $v_h= P_h v$ and $\tau_h=I_h\sigma$ in (\ref{theorem-err1}),
  by (\ref{proj-error}) and \meqref{p-err2}, we
   obtain  \meqref{t1}.
\end{proof}

\begin{remark}  It immediately follows from  Theorem \ref{equal} and Lemma \ref{lemma1} that there exists an interpolation operator
 $\Pi_h: H^1(\Omega, \mathbb{S})\rightarrow \Sigma_{k,h}$ such that
$$
 (\d(\tau-\Pi_h\tau), v_h)_K=0\text{  for  any } K \text{ and }v_h\in V_{k,h}
$$
 for any $\tau\in H^1(\Omega, \mathbb{S})$.  Further, if $\tau\in H^{k+1}(\Omega, \mathbb{S})$, it holds that
 $$
 \|\tau-\Pi_h\tau\|_{L^2(\Omega)}\leq Ch^{k+1}|\tau|_{H^{k+1}(\Omega)}.
 $$
 This, together with a standard argument,  leads to the following optimal error estimate for the stress in the $L^2$ norm
 $$
 \|\sigma-\sigma_h\|_{L^2(\Omega)}\leq Ch^{k+1}|\sigma|_{H^{k+1}(\Omega)},
 $$
 provided that $\sigma\in H^{k+1}(\Omega, \mathbb{S})$.
\end{remark}

\begin{remark}
The extension to nearly incompressible or incompressible elastic materials is possible. In the homogeneous
isotropic case the compliance tensor is given by
\a{
      A \tau &= \frac 1{2\mu} \left(
       \tau - \frac{\lambda}{2\mu + n \lambda} \operatorname{tr}(\tau)
        \delta \right), }
  where $\delta$ is an identity matrix, and  $\mu > 0$, $\lambda>0$ are the Lam\'{e} constants. For our mixed
method, as for most methods based on the Hellinger--Reissner principle,
one can prove that the error estimates hold uniformly in $\lambda$.
In the analysis above we use the fact that
$$
\alpha\|\tau\|_0\leq (A\tau, \tau)
$$
for some positive constant $\alpha$. This estimate degenerates
$\alpha\rightarrow 0$ when $\lambda\rightarrow +\infty$.
However the estimate remains true with $\alpha>0$  depending only on $\Omega$
 and $\mu$ if we restrict $\tau$ to functions for which $\d \tau=0$ and
 $\int_{\Omega}\text{tr}(\tau) d\b x=0$, see \cite{Brezzi-Fortin}, also \cite{Arnold-Douglas-Gupta,XieXu2011} for more details.
 \end{remark}

\section{Conclusions}
In this paper we propose a family of mixed elements of  symmetric tensors in any dimension.
 For stability, we require  in  Section \ref{sec3} that the polynomial degree for the stress be greater than $n$. Note that one key result, namely,  Theorem \ref{equal} holds for an arbitrary $k$, which, in a forth coming paper, will be used  to design lower order methods such that $1\leq k\leq n$.  In addition, the results in this paper will be used, in that paper,   to  derive, in a unified way,  those elements in  \cite{Arnold-Winther-conforming} and \cite{Arnold-Awanou-Winther},  and generalize them to any space dimension.

\section*{Acknowledgement} The author would like to thank Professor Jinchao Xu  for his constructive suggestions, in particular,  his  valuable
 suggestion that the author  define a set of local  degrees of freedom of shape function spaces for stress on each element.


\begin{thebibliography}{999}

\bibitem{Adams} R. A. Adams, Sobolev Spaces, New York: Academic
Press, 1975.

\bibitem{Adams-Cockburn} S. Adams and B. Cockburn,
  A mixed finite element method for elasticity in three dimensions, J.
    Sci. Comput. 25 (2005),  515--521.

\bibitem{Amara-Thomas} M. Amara and J. M. Thomas,
   Equilibrium finite elements for the linear elastic problem, Numer.
    Math. 33 (1979), 367--383.

\bibitem{Arnold2002} D. N. Arnold,  Proceedings of the International Congress of Mathematicians, Vol. I: Plenary Lectures and Ceremonies (2002), 137-157.

\bibitem{Arnold-Awanou} D. N. Arnold and G. Awanou,
   Rectangular mixed finite elements for elasticity, Math. Models
    Methods Appl. Sci. 15 (2005),  1417--1429.

\bibitem{Arnold-Awanou-Winther}
    D. Arnold, G. Awanou and R. Winther,
   Finite elements for symmetric tensors in three dimensions,
    Math. Comp. 77 (2008),  1229--1251.

\bibitem{Arnold-Brezzi-Douglas}
     D. N. Arnold, F. Brezzi and J. Douglas, Jr.,
     PEERS: A new mixed finite element for plane elasticity,
    Jpn. J. Appl. Math. 1 (1984),  347--367.

\bibitem{Arnold-Douglas-Gupta}
  D. N. Arnold, J. Douglas Jr., and C. P. Gupta,
    A family of higher order mixed finite element
    methods for plane elasticity, Numer. Math. 45 (1984),  1--22.

\bibitem{Arnold-Falk-Winther}
  D. N. Arnold, R. Falk and R. Winther,
  Mixed finite element methods for linear elasticity with
    weakly imposed symmetry, Math. Comp. 76 (2007),  1699--1723.


\bibitem{Arnold-Winther-conforming}
    D. N. Arnold and R. Winther,
    Mixed finite element for elasticity, Numer. Math. 92 (2002),
     401--419.

\bibitem{Arnold-Winther-n} D. N. Arnold and R. Winther,
   Nonconforming mixed elements for elasticity, Math. Models.
   Methods Appl. Sci. 13 (2003), 295--307.


\bibitem{Awanou} G. Awanou, Two remarks on rectangular mixed finite elements for elasticity, J. Sci.
 Comput. 50 (2012), 91--102.


\bibitem{Boffi-Brezzi-Fortin}
   D. Boffi, F. Brezzi and M. Fortin,
   Reduced symmetry elements in linear elasticity, Commun.
    Pure Appl. Anal. 8 (2009),   95--121.
\bibitem{Brezzi}
   F. Brezzi, On the existence, uniqueness and approximation of
   saddle-point problems arising from Lagrangian multipliers, Rev.
   Francaise Automat. Informat. Recherche Operationnelle Ser.
   Rouge, 8(R-2) (1974), 129--151.

\bibitem{Brezzi-Fortin}   F. Brezzi and  M. Fortin,
   Mixed and hybrid finite element methods, Springer, 1991.

\bibitem{Carstensen-Eigel-Gedicke}C. Carstensen, M. Eigel, J. Gedicke,
  Computational competition of symmetric mixed FEM in linear
  elasticity, Comput. Methods Appl. Mech. Engrg. 200 (2011),
   2903--2915.

\bibitem{Carstensen-Gunther-Reininghaus-Thiele2008}
  C. Carstensen,  D. G\"{u}nther, J. Reininghaus, J. Thiele,
  The Arnold--Winther mixed FEM in linear elasticity. Part I:
  Implementation and numerical verification, Comput. Methods Appl.
   Mech. Engrg. 197 (2008), 3014--3023.

\bibitem{Chen-Wang}   S. C. Chen and Y. N. Wang,
Conforming rectangular mixed finite elements for elasticity, J.
Sci. Comput. 47 (2011),  93--108.

%\bibitem{Clement}   P. Cl\'ement,
%   Approximation by finite element functions using local regularization,
%   RAIRO Anal. Numer., 2 (1975), 77--84.

\bibitem{Duran2001}R. G. Dur\'{a}n and M. A. Muschietti, An explicit right inverse of the divergence
operator which is continuous in weighted norms, Studia Mathematica,148(2001), 207--219.


%http://mate.dm.uba.ar/~rduran/papers/dm.pdf

\bibitem{Cockburn-Gopalakrishnan-Guzman}
  B. Cockburn, J. Gopalakrishnan and J. Guzm\'an,
    A new elasticity element made for enforcing weak stress symmetry,
    Math. Comp. 79 (2010),  1331--1349.

%\bibitem{Fraejis}
% B. M. Fraejis de Veubeke, Displacement and equilibrium models in
% the finite element method, Stress analysis, (O. C. Zienkiewics and G. S. Holister, eds.)
% New York, Wiley (1965), PP. 145--197.

\bibitem{GiraultRaviart1986} V. Girault, P. A. Raviart, Finite Element Methods for Navier-Stokes equations,
Springer, Berlin, 1986.

\bibitem{Gopalakrishnan-Guzman-n} J. Gopalakrishnan and J. Guzm\'an,
  Symmetric nonconforming mixed finite elements for linear elasticity,
 SIAM J. Numer. Anal. 49 (2011),  1504--1520.

\bibitem{Gopalakrishnan-Guzman} J. Gopalakrishnan and J. Guzm\'an,
    A second elasticity element using the matrix bubble,
    IMA J. Numer. Anal. 32 (2012),  352--372.


\bibitem{Guzman} J. Guzm\'an,
   A unified analysis of several mixed methods for elasticity
   with weak stress symmetry,
    J. Sci. Comput. 44 (2010),  156--169.

\bibitem{Hu2013}J. Hu,  A new family of efficient  conforming mixed
elements on both rectangular and cuboid meshes for linear  elasticity in the symmetric formulation,
arXiv:1311.4718v3 [math.NA], 17 Dec 2013

\bibitem{Hu-Shi}  J. Hu and Z. C. Shi,
     Lower order rectangular nonconforming mixed elements for plane elasticity,
     SIAM J. Numer. Anal.  46 (2007),   88--102.

\bibitem{Hu-Man-Zhang2013} J. Hu, H. Y. Man and S. Zhang,
    The minimal mixed finite element method for the symmetric stress
   field on rectangular  grids in any space dimension, arXiv:1304.5428[math.NA] (2013).

\bibitem{Hu-Man-Zhang2014}J. Hu, H. Y. Man and S. Zhang,
A simple conforming mixed finite element for linear
elasticity on rectangular grids in any space dimension,
   J. Sci. Comput. 58(2014), 367--379.

\bibitem{Hu-Zhang2014a} J. Hu and S. Zhang, A family of conforming mixed finite elements for linear elasticity on
 triangle grids, arXiv:1406.7457 [math.NA], 2014.

 \bibitem{Hu-Zhang2014b} J. Hu and S. Zhang, A family of conforming mixed finite elements for linear elasticity on
  tetrahedral grids, Sci. China  Math.,58(2015), pp. 297--307; see also arXiv:1407.4190 [math.NA], 2014.


\bibitem{Johnson-Mercier} C. Johnson and B. Mercier,
    Some equilibrium finite element methods for two-dimensional elasticity
     problems, Numer.Math. 30 (1978),  103--116.

\bibitem{Man-Hu-Shi} H.-Y. Man, J. Hu and Z.-C. Shi,
   Lower order rectangular nonconforming mixed finite element for
     the three-dimensional elasticity problem,
   Math. Models Methods Appl. Sci. 19 (2009),  51--65.

\bibitem{Morley} M. Morley, A family of mixed finite elements
   for linear elasticity, Numer. Math. 55 (1989),  633--666.

\bibitem{Scott-Zhang}
   L. R. Scott and S. Zhang,
     Finite-element interpolation of
   non-smooth functions satisfying boundary conditions, Math. Comp.
     54  (1990), 483--493.

\bibitem{Stenberg-1}
     R. Stenberg, On the construction of optimal mixed finite
     element methods for the linear elasticity
     problem, Numer. Math. 48 (1986),   447--462.

\bibitem{Stenberg-2}
     R. Stenberg, Two low-order mixed methods for the elasticity problem,
     In: J. R. Whiteman (ed.):
     The Mathematics of Finite Elements and Applications, VI. London:
     Academic Press, 1988, 271--280.
\bibitem{Stenberg-3} R. Stenberg,
   A family of mixed finite elements for the elasticity problem,
    Numer. Math. 53 (1988),   513--538.


%\bibitem{Watwood}
% V. B. Watwood Jr. and B. J. Hartz, An equilibrium stress field
% model for finite element solution of two-dimensional elastostatic
% problems, Internat. Jour. Solids and Structures, 4 (1968),
% 857--873.
\bibitem{XieXu2011}X. P.  Xie and J. C. Xu, New mixed finite elements for plane elasticity
and Stokes equations, Sci. China Math.,54(2011), 1499--1519.



\bibitem{Yi-3D} S. Y. Yi,
    Nonconforming mixed finite element methods for
       linear elasticity using rectangular
     elements in two and three dimensions,
      CALCOLO 42 (2005), 115--133.

\bibitem{Yi} S. Y. Yi,
    A New nonconforming mixed finite element method for linear elasticity,
    Math. Models.Methods Appl. Sci. 16 (2006),  979--999.

\bibitem{ZTZ2005}O. C. Zienkiewicz, R. L. Taylor, and  J. Z. Zhu, The Finite Element Method: Its  Basis and
         Fundamentals, 6th ed., vol. 1, Amsterdam--Boston--Heidelberg--London--New York--Oxford--Paris--San Diego--San Francisco--Singapore--Sydney--Tokyo, 2005.

\end{thebibliography}
\end{document}